\documentclass[12pt]{amsart}

\usepackage{graphicx,url,lineno}
\usepackage{amssymb}
\usepackage{epstopdf}
\usepackage{color}
\usepackage{mathrsfs}

\newtheorem{lemma}{Lemma}
\newtheorem{definition}{Definition}
\newtheorem{remark}{Remark}

\newtheorem{proposition}{Proposition}

\title[A note on the fractional logistic equation]{A note on the fractional logistic equation}

\author[Area]{I. Area}
\author[Losada]{J. Losada}
\author[Nieto]{J.\,J. Nieto}

\address[Area]{Departamento de Matem\'{a}tica Aplicada II, E.E. Telecomunicaci\'{o}n, Universidade de Vigo, 36310-Vigo, Spain.}\email[Area]{area@uvigo.es}

\address[Losada]{Facultade de Matem\'{a}ticas, Universidade de Santiago de
Compostela, 15782-Santiago de Compostela, Spain. Corresponding author}
\email[Losada]{jorge.losada@rai.usc.es}

\address[Nieto]{Facultade de Matem\'{a}ticas, Universidade de Santiago de
Compostela, 15782-Santiago de Compostela, Spain, and Faculty of Science, King Abdulaziz University, P.O. Box 80203, 21589, Jeddah, Saudi Arabia.}
\email[Nieto]{juanjose.nieto.roig@usc.es}

\date{\today}                                        

\begin{document}

\begin{abstract}
In this short note, we show that the real function recently proposed by Bruce J. West [Exact solution to fractional logistic equation, Physica A 429 (2015) 103--108] is not an exact solution for the fractional logistic equation.
\end{abstract}

\maketitle

\section{Introduction}\label{S:1}

The topic of Fractional Calculus (that is, the study of integrals and derivatives of non-integer order) is nowadays a growing interest area of mathematics. Moreover, physicists and engineers are also really interested in applications coming from this nice theory which origin goes back to the theory of differential calculus; more precisely to Leibniz's famous note in his letter to L'H\^opital dated 30 September 1695.

For a long time, the theory of Fractional Calculus developed only as a theoretical field of mathematics. However, in the last decades, it was shown that some fractional operators describe in a better way some complex physical phenomena, especially when dealing with memory processes or viscoelastic and viscoplastic materials. Well known references about the application of fractional operators in rheology modeling are \cite{bagley1,bagley2}. One of the most important advantage of fractional order models in comparison with integer order ones is that fractional integrals and derivatives are a powerful tool for the description of memory and hereditary properties of some materials. Notice that integer order derivatives are local operators, but the fractional order derivative of a function in a point depends on the past values of such function. This features motivated the successful use of fractional calculus in CRONE \cite{oustaloup} and PID controllers \cite[Chapter 9]{podlubny}.

Due to this, in the last two or three decades, a great interest has been devoted to the study of fractional differential equations. Thus, fractional order differential equations (one day called extraordinary differential equations) play today a very important role describing some real world phenomena. For a complete exposition of the theory of Fractional Calculus one can see \cite{trujillo,spanier,podlubny,samko}.\medbreak

The exponential function, $\exp(t)$, plays a fundamental role in mathematics and it is really useful in the theory of integer order differential equations. In the case of fractional order, it loses some beautiful properties and Mittag-Leffler function appears as its natural substitute. Next, we recall its definition and some basic properties.
\begin{definition}
The function $E_\alpha(z)$ is named after the great Swedish mathematician G\"osta Mittag-Leffler (1846-1927) who defined it as a power series given by
\[E_\alpha(z)=\sum_{k=0}^\infty\dfrac{z^k}{\Gamma(\alpha k+1)},\qquad\alpha>0.\]
\end{definition}
This function provides a simple generalization of the exponential function because of the replacement of $k!=\Gamma(k+1)$ by $(\alpha k)!=\Gamma(\alpha k+1)$ in the denominator of the terms of the exponential series. Due to this, such function can be considered the simplest nontrivial generalization of exponential function.

\begin{definition}
A two parameter function of Mittag-Leffler type is defined by
\[E_{\alpha,\beta}(z)=\sum_{k=0}^\infty\dfrac{z^k}{\Gamma(\alpha k+\beta)},\qquad \alpha>0,\,\beta>0.\]
\end{definition}

During the first half of the twntieth century, the Mittag-Leffler functions remained almost unkonown to the majority of the mathematicians. However, recently attention of mathematicians and other scientists towards functions of Mittag-Leffler type has increased and today some mathematicians like to refer the classical Mittag-Leffler function as the {\it{Queen Function}} of Fractional Calculus, and to consider all related functions as her court.

Recently, a complete monograph (see \cite{gorenflo}) about Mittag-Leffler functions and its several applications has been published. More information about this kind of special functions ans its relations with Fractional Calculus can be found in \cite{mainardi}.

It follows from previous definitions that, for instance,

\begin{align*}
&E_{1,1}(z)=\exp(z),&E_{2,1}(z^2)=\cosh(z),\\
&E_{1,2}(z)=\dfrac{\exp(z)-1}{z},&E_{2,2}(z^2)=\dfrac{\sinh(z)}{z},\\
&E_{1,3}(z)=\dfrac{\exp(z)-1-z}{z^2},&E_{2,2}(z^2)=\cos(z),
\end{align*}
which proves the importance of Mittag-Leffler functions in mathematics. Furthermore, using the following well known result
\[\left[{}^{\textsc{c}}D^\alpha x^b\right]\,(t)=t^{b-\alpha}\dfrac{\Gamma(b+1)}{\Gamma(b+1-\alpha)},\qquad b>-1,\,b\neq 0,\]
we can easily see that $E_\alpha(\lambda t^\alpha)$, with $\lambda\in\mathbb{R}$, is a nontrivial eigenfunction for Caputo fractional derivative operator (previously denoted by $^{\textsc{c}}D^\alpha$, see \cite[Section 2.4]{trujillo}). That is, if $f(t)=E_\alpha(\lambda t^\alpha)$ for $t\ge 0$ then
\[^{\textsc{c}}D^\alpha f\,(t)=\lambda f(t),\qquad t\ge 0.\]\medbreak

In the past, some authors used to derived some results the following {\it{property}} of Mittag-Leffler function
\[E_\alpha(a(t+s)^\alpha)=E_\alpha(at^\alpha)E_\alpha(as^\alpha),\]
where $a$ is a real constant. Unfortunately, such property is unavailable unless $\alpha=1$ or $a=0$, both of them are trivial situations, since we obtain exponential function in the first case, and a constant function secondly. For more information about the invalidity of such equality we refer to \cite{peng}.

In this note we discuss a similar property related again with Mittag-Leffler function and using such property, we conclude that the function introduced in \cite{west} is not an exact solution for the fractional logistic equation. It may be obvious, but nevertheless Mittag-Leffler function shares some analogue properties with exponential function, these kind of results show that exponential function cannot be replaced, in general, by Mittag-Leffler function in the theory of Fractional Calculus.

\section{Main results}\label{S:2}

In 1838, P.~F.~Verhulst introduced a nonlinear term into the rate equation; he was studying population models and he wanted to avoid the catastrophic predictions previously proposed by T.~Malthus, who had used the rate equation to model human population growth (for historical references see \cite{west}). By this way, P.~F.~Verhulst obtained what today is known as the logistic equation:
\begin{equation}\label{ec:1}
u'(t)=ku(t)(1-u(t)),\qquad t\ge 0.
\end{equation}
This differential equation is one of the few nonlinear differential equations that has a known exact closed form solution, which is given by
\begin{equation}\label{ec:2}
u(t)=\dfrac{u_0}{u_0+(1-u_0)\exp(-kt)},\qquad t\ge 0,
\end{equation}
where $u_0$ is the initial state (that is, $u(0)=u_0=N(0)/N_{\rm max}$, where $N(0)$ is the total population at the initial time and $N_{\rm max}$ is the carrying capacity of the ecosystem). Since that time, the logistic equation has found several applications ranging from machine learning to modeling growth of tumors.\medbreak

Fractional logistic equation has been studied, see for instance \cite{discrete,sayed}; but until now no exact is known. Recently, Bruce J. West has published a research paper (see \cite{west}) where, using the Carleman embeding technique, he proposes an exact solution for the fractional logistic equation. More concretely, in such reference it is said that the function given by
\begin{equation}\label{ec:3}
u(t)=\sum_{n=0}^\infty\left(\dfrac{u_0-1}{u_0}\right)E_\alpha(-nk^\alpha t^\alpha),\qquad t\ge 0,
\end{equation}
may be an exact solution of the fractional differential equation 
\begin{equation}\label{ec:4}
^{\textsc c}D^\alpha u\,(t)=k^\alpha u(t)(1-u(t)),\qquad t\ge 0,\qquad 0<\alpha<1.
\end{equation}

Of course, following \cite{west}, notice that if $\alpha=1$ then function given in (\ref{ec:3}) is equivalent to
\[u(t)=\sum_{n=0}^\infty\left(\dfrac{u_0-1}{u_0}\right)\exp(-nk t)=\dfrac{u_0}{u_0+(1-u_0)\exp(-kt)},\qquad t\ge 0\]
which coincides with (\ref{ec:2}). That is, if in (\ref{ec:3}) $\alpha=1$, then we obtain the solution of the integer order logistic equation (\ref{ec:1}). Moreover, fractional solutions given by (\ref{ec:3}) have similar properties to integer order solution given by (\ref{ec:2}).

\begin{remark}
Since Carleman embedding technique was introduced to solve integer order differential equations, not fractional order differential equations, the author of \cite{west} warms us about the necessity to validate solution of (\ref{ec:4}) given by (\ref{ec:3}). This is the main objective of this short note.
\end{remark}\medbreak

In what follows, we assume that function given by (\ref{ec:3}) is an exact solution of the fractional differential equation (\ref{ec:4}). 

\begin{lemma}\label{L:2.1}
If the function given by (\ref{ec:3}) is an exact solution of fractional logistic equation (\ref{ec:4}) of order $0<\alpha<1$, then for all $n\in\mathbb{N}\cup\{0\}$ we have that
\begin{equation}\label{ec:5}
(n+1)E_\alpha(-nk^\alpha s^\alpha)=\sum_{j=0}^nE_\alpha(-(n-j)k^\alpha t^\alpha)E_\alpha(-jk^\alpha t^\alpha).
\end{equation}
\begin{proof}
Consider for each $j\in \mathbb{N}\cup\{0\}$
\[a_j=\left(\dfrac{u_0-1}{u_0}\right)^jE_\alpha(-jk^\alpha t^\alpha).\]
Thus, function given by (\ref{ec:3}) can be written as
\[u(t)=\sum_{n=0}^\infty a_n,\qquad t\ge 0.\]
Moreover, given $0<\alpha<1$, we have that
\[^{\textsc c}D^\alpha\left[\sum_{n=0}^\infty\left(\dfrac{u_0-1}{u_0}\right)^n E_\alpha(-nk^\alpha s^\alpha)\right]\,(t)=-k^\alpha\sum_{n=0}^\infty\left(\dfrac{u_0-1}{u_0}\right)^nn E_\alpha(-nk^\alpha s^\alpha),\]
and using the Cauchy product, we obtain:

\begin{align*}
u(t)(1-u(t))&=\sum_{n=0}^\infty\left(\dfrac{u_0-1}{u_0}\right)^nE_\alpha(-nk^\alpha t^\alpha)\left(1-\sum_{n=0}^\infty\left(\dfrac{u_0-1}{u_0}\right)^nE_\alpha(-nk^\alpha t^\alpha)\right)\\
&=\sum_{n=0}^\infty\left(\dfrac{u_0-1}{u_0}\right)^nE_\alpha(-nk^\alpha t^\alpha)-\sum_{n=0}^\infty\sum_{j=0}^n a_ja_{n-j}\\
&=\sum_{n=0}^\infty\left(\dfrac{u_0-1}{u_0}\right)^n\left[E_\alpha(-nk^\alpha t^\alpha)-\sum_{j=0}^nE_\alpha(-(n-j)k^\alpha t^\alpha)E_\alpha(-jk^\alpha t^\alpha)\right].
\end{align*}
Since $u(t)$ is a solution of equation (\ref{ec:4}), we deduce that for all $n\in\mathbb{N}\cup\{0\}$ it must be
\begin{equation*}
(n+1)E_\alpha(-nk^\alpha s^\alpha)=\sum_{j=0}^nE_\alpha(-(n-j)k^\alpha t^\alpha)E_\alpha(-jk^\alpha t^\alpha),
\end{equation*}
which is the equality given in the statement of the lemma.\qed
\end{proof}
\end{lemma}

\begin{proposition}\label{P:2.2}
In the condition of Lemma \ref{L:2.1}, we have that
\begin{equation}\label{ec:6}
E_\alpha\left(-2k^\alpha t^\alpha\right)=E_\alpha\left(-k^\alpha t^\alpha\right)E_\alpha\left(-k^\alpha t^\alpha\right).
\end{equation}
\begin{proof}
Considering $n=2$ in equality (\ref{ec:5}), we obtain the desired result.\qed
\end{proof}
\end{proposition}

\begin{proposition}\label{P:2.3}
Identity (\ref{ec:6}) of Proposition \ref{P:2.2} holds if and only if $\alpha=1$.
\begin{proof}
By definition of the Mittag-Leffler function we have that,
\begin{equation}\label{ec:7}
E_\alpha(-2k^\alpha t^\alpha)=\sum_{n=0}^\infty\dfrac{(-2k^\alpha t^\alpha)^n}{\Gamma(n\alpha+1)},\qquad t\ge 0.
\end{equation}
Moreover, using the Cauchy product, we obtain:
\begin{align}\label{ec:8}
\nonumber
E_{\alpha}(-k^\alpha t^\alpha)E_\alpha(-k^\alpha t^\alpha)&=\sum_{n=0}^\infty\sum_{j=0}^n\dfrac{(-k^\alpha t^\alpha)^{n-j}(-k^\alpha t^\alpha)^j}{\Gamma((n-j)\alpha+1)\Gamma(j\alpha+1)}\\
\nonumber
&=\sum_{n=0}^\infty(-k^\alpha t^\alpha)^n\sum_{j=0}^n\dfrac{1}{\Gamma((n-j)\alpha+1)\Gamma(j\alpha+1)}\\
\nonumber
&=\sum_{n=0}^\infty(-2k^\alpha t^\alpha)^n\sum_{j=0}^n\dfrac{2^{-n}}{\Gamma((n-j)\alpha+1)\Gamma(j\alpha+1)}\\
&=\sum_{n=0}^\infty(-2k^\alpha t^\alpha)^nb_n,\qquad t\ge 0,
\end{align}
where
\[b_n=\sum_{j=0}^n\dfrac{2^{-n}}{\Gamma((n-j)\alpha+1)\Gamma(j\alpha+1)}.\]
Equality (\ref{ec:6}) implies that coefficents of series expansions (\ref{ec:7}) and (\ref{ec:8}) must be equal. Particularly, if we choose $n=2$, we must have
\begin{align*}
a_2=\dfrac{1}{4}\left(\dfrac{2}{\Gamma(2\alpha+1)}+\dfrac{1}{\Gamma(\alpha+1)\Gamma(\alpha+1)}\right)=\dfrac{1}{\Gamma(2\alpha+1)},
\end{align*}
or equivalently,
\begin{equation}\label{ec:9}
\dfrac{\Gamma(2\alpha+1)}{4\Gamma(\alpha+1)\Gamma(\alpha+1)}=\dfrac{1}{2}.
\end{equation}
Since if $\alpha\in(0,1)$ then $\Gamma(2\alpha+1)<\Gamma(3)=2$ and $4\Gamma(\alpha+1)^2>4\Gamma(1)^2=4$, it is clear that equality in (\ref{ec:9}) holds if and only if $\alpha=1$.\qed
\end{proof}
\end{proposition}

\begin{remark}
From Proposition \ref{P:2.3}, it is clear that function given by (\ref{ec:3}) is a solution of differential equation (\ref{ec:4}) if and only if $\alpha=1$, but in such case, we have an integer order differential equation and function defined by (\ref{ec:3}) coincides with the solution of such differential equation.

Notice also that since when $\alpha=1$ then $E_\alpha(t)=\exp(t)$, identities (\ref{ec:5}) and (\ref{ec:6}) are well known properties.
\end{remark}

\begin{remark}
Differentiating in both sides of equation (\ref{ec:6}), we deduce another identity which only holds if $\alpha=1$,
\[E_{\alpha,\alpha}\left(-2k^\alpha t^\alpha\right)=E_{\alpha,\alpha}\left(-k^\alpha t^\alpha\right)E_\alpha\left(-k^\alpha t^\alpha\right).\]
\end{remark}\medbreak

Using Mathematica \cite{mathematica} (see also \cite{computation}), we can reaffirm analytical results previously obtained. In figure 1 some graphs can be found in this direction. Notice that when $\alpha$ is near $1$ the difference is smaller than in the case $\alpha$ near $0$.\medskip

\begin{figure}[h!]
\begin{center}
\includegraphics[width=0.44\textwidth]{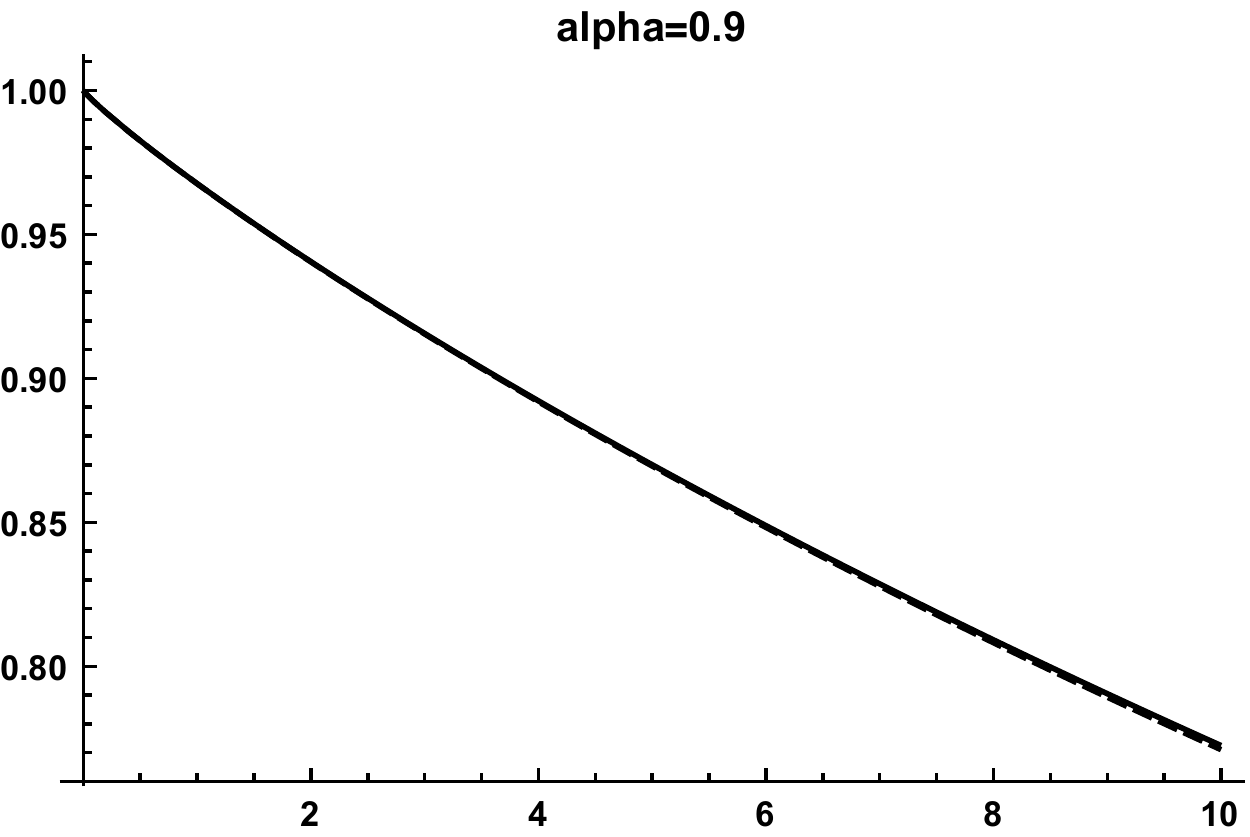}\hspace*{0.4cm}
\includegraphics[width=0.44\textwidth]{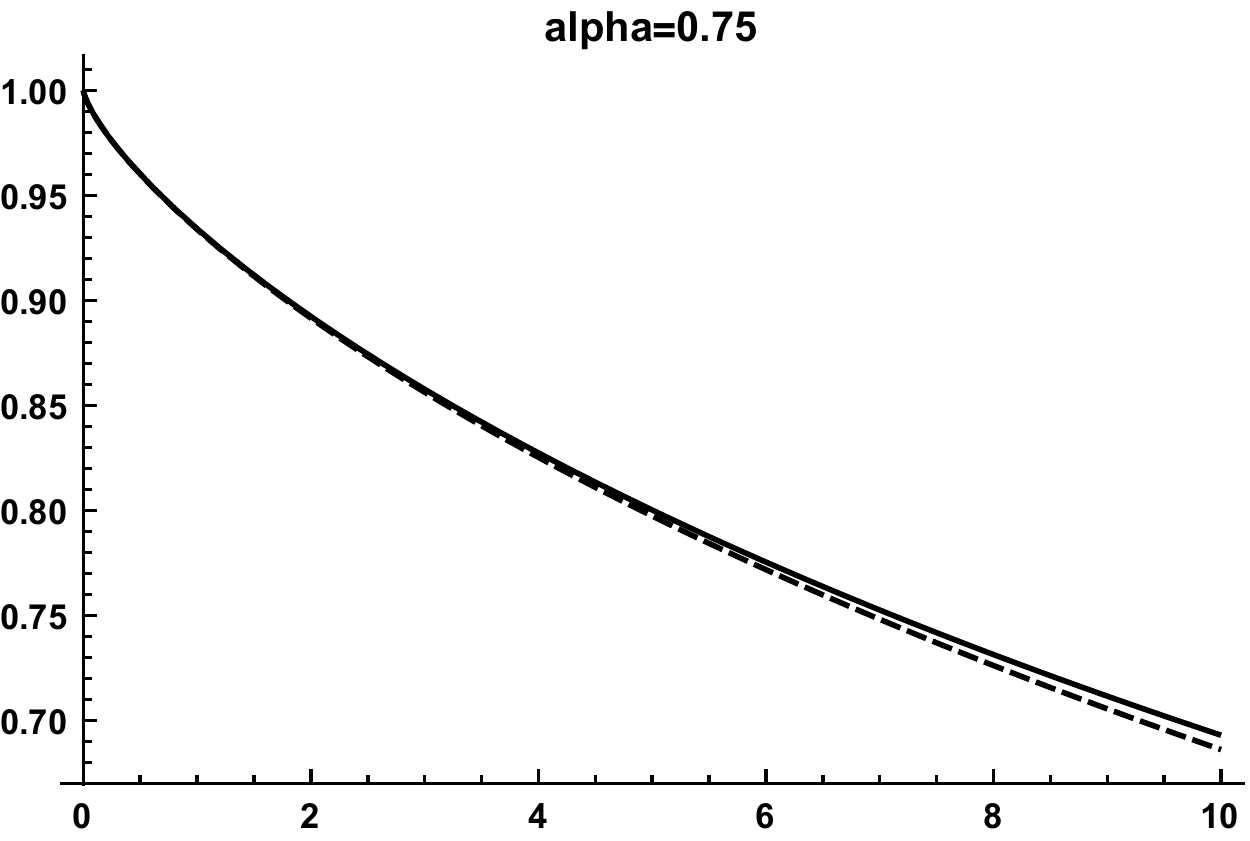}\\
\includegraphics[width=0.44\textwidth]{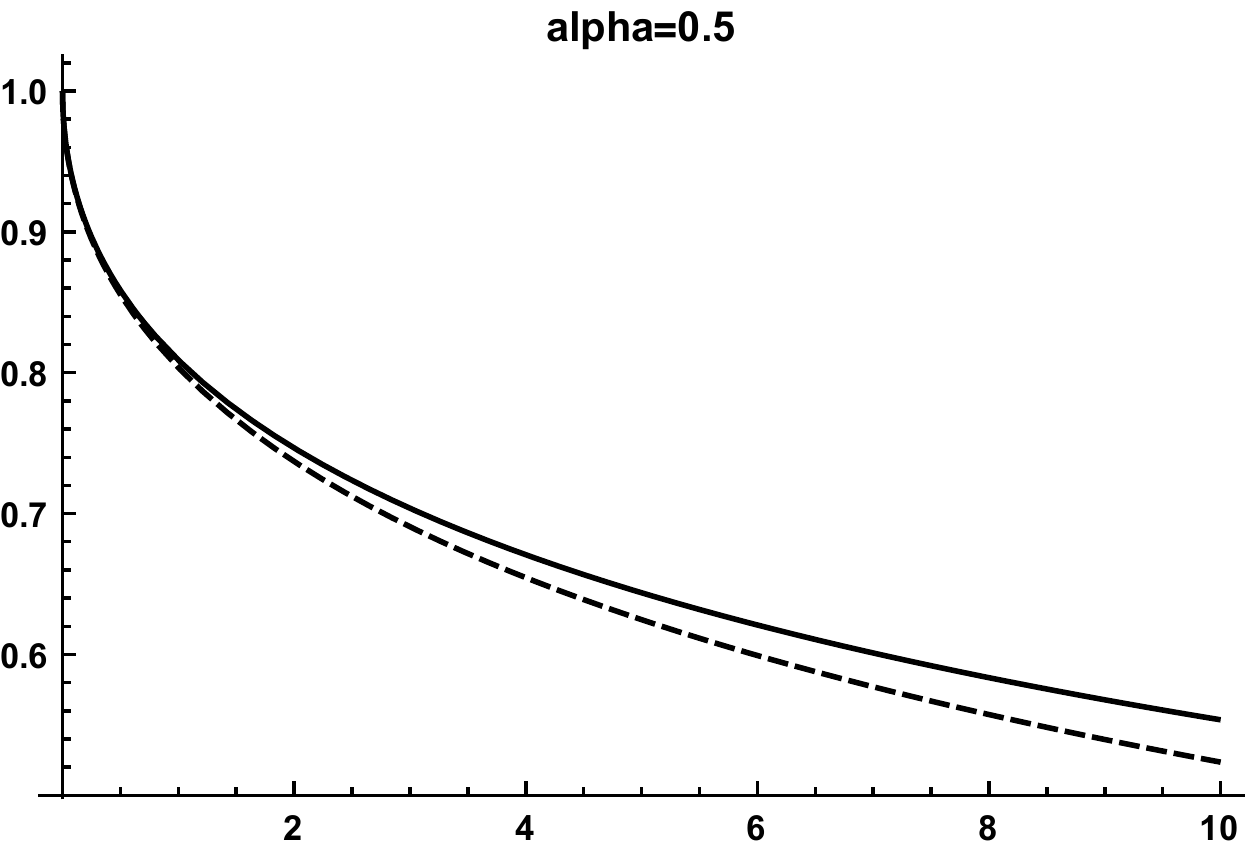}\hspace*{0.4cm}
\includegraphics[width=0.44\textwidth]{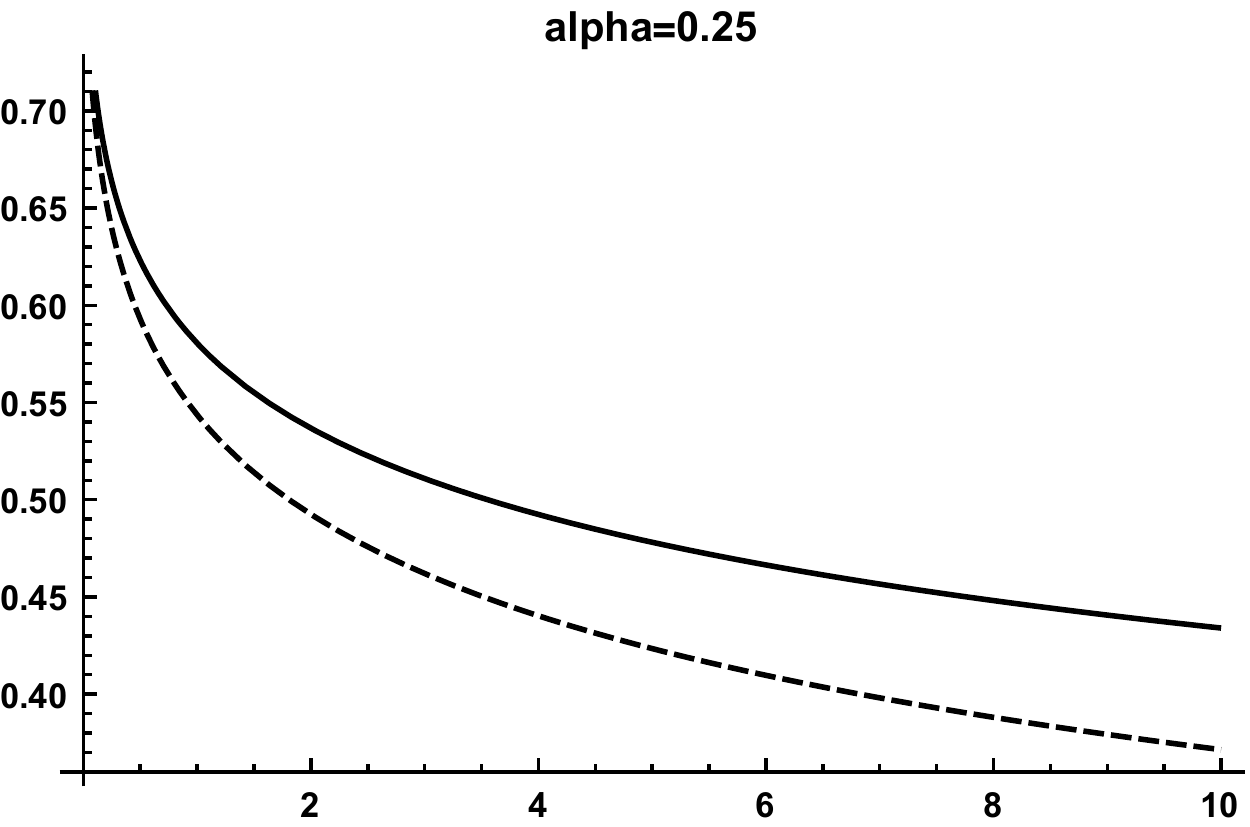}
\caption{Comparison considering diferent values of $\alpha$ ($\alpha=0.9,0.75,0.5,0.25$) between $E_\alpha\left(-2k^\alpha t^\alpha\right)$ (continuous line) and the product $E_\alpha\left(-k^\alpha t^\alpha\right)E_\alpha\left(-k^\alpha t^\alpha\right)$ (dashed line).}
\end{center}
\end{figure}

Recently, Turalska and West \cite{turalska} have used a spectral technique to study the solution to the fractional logistic differential equation since the exact solution is not yet known. For the approximate siolution there is a small, apparently generic, systematic error that they are not able to fully interpret.

\section{Conclusion}
In this short note, we have  shown that the function recently proposed by Bruce J. West is not an exact solution for the fractional logistic equation.

In any case, we would like to mention that the use of Carleman embeding technique may be a powerful tool in the study of some fractional differential equations, but its application must be studied carefully.

\section*{Acknowledgments}
We would like to express our sincere gratitude and appreciation to Bruce J. West for his helpful and valuable comments concurring with our findings.\medskip

The work of I. Area has been partially supported by the Ministerio de Econom\'{\i}a y Competi\-tividad of Spain under grant MTM2012--38794--C02--01, co-financed by the European Community fund FEDER. J.J. Nieto also acknowledges partial financial support by the Ministerio de Econom\'{\i}a y Competitividad of Spain under grant MTM2013--43014--P co-financed by the European Community fund FEDER. J. Losada  acknowledges financial support by Xunta de Galicia under grant Plan I2C ED481A-2015/272.

\end{document}